\pdfoutput=1
\documentclass[a4paper, 11pt, abstracton, final]{scrartcl}
\usepackage[utf8]{inputenc}
\usepackage{cite}
\usepackage{amsmath}
\usepackage{amssymb}
\usepackage{listings}
\usepackage{algorithm}
\usepackage[bookmarks=false]{hyperref}
\usepackage{enumitem} 

\usepackage[utf8]{inputenc}
\usepackage{tikz}
\usetikzlibrary{arrows, decorations.pathreplacing}

\usepackage{environ}
\usepackage[textsize=footnotesize,textwidth=5cm]{todonotes}

\newcommand\todoblock[3][]{\todo[inline,prepend,caption={\textbf{#2}},#1]{#3}}

\NewEnviron{todoenv}[2][]{\todoblock[#1]{#2}{\BODY}}
\NewEnviron{todolist}{\begin{todoenv}\begin{itemize} \BODY \end{itemize}\end{todoenv}}
\newtheorem{example}{Example}

\usepackage{xcolor}

\hypersetup{
    colorlinks,
    linkcolor={red!50!black},
    citecolor={blue!50!black},
    urlcolor={blue!80!black}
}

\newcommand{\RR}{\mathbb{R}}
\newcommand{\NN}{\mathbb{N}}
\newcommand{\sgn}{\operatorname{sgn}}
\newcommand{\M}{\mathcal{M}}
\newcommand{\cell}{E}
\newcommand{\face}{\gamma}
\newcommand{\edge}{\kappa}
\newcommand{\vertex}{\xi}

\newcommand{\MC}{\textsf{MC}}
\newcommand{\MCxx}{\textsf{MC33}}
\newcommand{\KMT}{\textsf{K/MT}}
\newcommand{\ourMC}{\textsf{TPMC}}
\newcommand{\Tr}{T_{\operatorname{rel}}}

\def\clap#1{\hbox to 0pt{\hss#1\hss}}

\begin{document}

\lstset{
basicstyle=\fontfamily{zlmtt}\selectfont, language=C++,
numbers=left,
numbersep=5pt,
numberstyle=\sffamily\tiny\color{gray},
commentstyle=\color{gray}\emph,
keywordstyle={\fontseries{b}\selectfont\color{gray}},
keywords={std, tpmc},
keywordstyle=[2]{\fontseries{b}\selectfont},
keywords=[2]{vector, array, Coordinate, static, unsigned, int, double,
typedef, size_t},
otherkeywords={::},
}

\title{Geometric Integration Over Irregular Domains with
  topologic Guarantees}
  \author{Christian Engwer\thanks{
    Institute for Computational und Applied Mathematics Einsteinstra\ss e 62, 48149 M\"unster, Germany;
  \{christian.engwer, andreas.nuessing\}@uni-muenster.de} \and Andreas Nüßing $^*$
}

\maketitle
\begin{abstract}
  Implicitly described domains are a well established tool in the
  simulation of time dependent problems, e.g. using level-set methods.
  In order to solve partial differential equations on such domains, a
  range of numerical methods was developed, e.g. the Immersed Boundary
  method, Unfitted Finite Element or Unfitted discontinuous Galerkin
  methods, eXtended or Generalised Finite Element methods, just to
  name a few. Many of these methods involve integration over cut-cells
  or their boundaries, as they are described by sub-domains of the
  original level-set mesh.

  We present a new algorithm to geometrically evaluate the integrals over
  domains described by a first-order, conforming level-set function.
  The integration is
  based on a polyhedral reconstruction of the implicit geometry,
  following the concepts of the Marching Cubes algorithm.
  The algorithm preserves various topological properties of the implicit geometry in
  its polyhedral reconstruction, making it suitable for Finite Element computations.
  Numerical experiments show second order accuracy of the
  integration.

  An implementation of the algorithm is available as free software,
  which allows for an easy incorporation into other projects. The
  software is in productive use within the DUNE framework~\cite{dune08:2}.

\end{abstract}

\section{Introduction}
\label{sec:introduction}

When dealing with partial differential equations (PDEs), methods using
implicitly described domains have established themselves as an
alternative to the traditional work flow using geometry adapted
meshes.  These methods usually solve the PDE on a larger domain and
the actual geometry information is incorporated later on.  Some
popular methods are the Immersed Boundary method \cite{peskin77},
Unfitted Finite Element \cite{Barrett:1987:FUF} or Unfitted
discontinuous Galerkin methods \cite{engwer2009udg}, eXtended
\cite{dolbow2001efe} or Generalised Finite Element methods
\cite{strouboulis2001gfe}.

One popular class of implicit domain descriptions is via a scalar function,
usually a so called level-set function \cite{osher1988levelset_first} or a phase-field
function \cite{fix1982phasefield_first}. These approaches are
particularly attractive, as they can easily describe moving domains
and even topology changes. In the following we have only very few
requirements on the implicit description. A domain $\Omega\subset \RR^d$ is described
using a continuous scalar function $\Phi: \hat\Omega \subset \RR^d \rightarrow \RR$ such that:
\begin{equation}
  \label{eq:levelset}
  \Phi(x)
  \begin{cases}
    < 0 &, \text{ if }x\in\Omega \\
    = 0 &, \text{ if }x\in\partial\Omega \\
    > 0 &, \text{ else.}
  \end{cases}
\end{equation}
In order to incorporate the geometry information into the
discretization, many of the aforementioned methods must integrate
a function over the domain, or the surface.
In order to evaluate these integrals, a discretized
representation of $\Omega$, or $\partial\Omega$, is needed.

An efficient algorithm to reconstruct the interface between $\Omega$ and its
complement is the \emph{marching cubes} algorithm presented in
\cite{Lorensen1987}. It originates from the field of computer graphics
and is used to visualize volume data on a structured hexahedral grid.
For each grid element a case number can be generated based on the sign
of $\Phi$ at the cell vertices. This key is than used as an index
into a look-up table where the triangulation of the interface for this
specific configuration is stored.
This algorithm can be applied straight forward to tetrahedra, known as
\emph{marching tetrahedra} and other geometries
\cite{elvins1992survey}. As the case for cubes is the most challenging
one, as we will discuss later, we only consider cubes throughout this
paper.

There are several cases where the interface can not be uniquely
determined based only on the sign of $\Phi$ at the corners. These
ambiguities can be resolved as presented in \cite{Chernyaev1995} with
the \emph{marching cubes 33} algorithm. For
each face being ambiguous it is checked if two diagonally
opposing vertices are connected over the face. In addition, a test
needs to be performed, if two vertices could be connected through the
element. For each case, the set of necessary tests is stored in a table. The
resulting triangulation of the interface is topologically
correct.
An implementation and completion of this algorithm has been published
in \cite{Lewiner2003}. As we will discuss later, also
  this implementation lacks support for some special cases in 3D,
  where the topology is not properly retained.

A method which uses this idea for geometric integration over irregular
domains has been presented in \cite{Min2007}. First cubes are
triangulated into simplices independent of $\Phi$. To these simplices,
a marching cubes type algorithm is applied to obtain a discrete
representation of $\Omega$ which is then used for integration.
Using only simplices has the benefit that the
amount of possible triangulation cases is significantly reduced. The
downside is that the multilinear level-set is projected onto a linear
level-set, resulting in a systematic distortion of the reconstruction
along the newly introduced diagonals of cubes. Even with a symmetric
criss-cross triangulation, these systematic error can not be avoided
completely.


In this paper we employ the marching cubes ideas for geometric
integration. We provide volume triangulations of the
interior and exterior domain, along with the interface triangulation.
These triangulations are topologically correct and consistent with
each other. This allows for geometric integration over implicitly
described domains.  The partitioning into interior and exterior of
each element's face based on the volume triangulation is consistent
with the triangulation which results from the application of the lower
dimensional algorithm to the face. Although we only discuss the case
of hexahedral elements, the algorithm is implemented for different
dimensions and element types.


In the next section we introduce mathematical
preliminaries and discuss the required consistencies. Section
\ref{sec:existing-approaches} elaborates the existing approaches, which
were already briefly mentioned in the introduction.
The major contribution of this paper is the presentation of a
\emph{topology
preserving marching cubes} (\ourMC) construction
in Section \ref{sec:algorithm},
where we also compare the
results of out new approach to the existing ones for chosen typical
pathologic cases. The method is implemented as a C++ library, using
python for automatic code generation as
discussed in
Section \ref{sec:implementation}.
Finally, in Section \ref{sec:evaluation} numerical
experiments are presented which verify the accuracy and the topologic consistency for the
\ourMC{} algorithm.

\section{Preliminaries}
\label{sec:preliminaries}

\subsection{Implicitly described domains}
\begin{figure}
  \centering
  \begin{tikzpicture}[scale=4]
  \def\radius{.25}
  \draw (.5,.5) node {$\Phi<0$};
  \draw (.25,.85) node {$\Phi>0$};
  \draw[gray] (.5,.5) ++(45:\radius) -- +(.075,.075) node[above, black] {$\Phi=0$};
  \draw (1,0.9) node[right] {\rlap{$\hat\Omega=[0,1]^2$}};
  \draw (0,0.9) node[left] {\rlap{~}};
  \draw[thick] (.5,.5) circle (\radius);
  \draw[thick] (0,0) rectangle (1,1);
\end{tikzpicture}
  \caption{Illustration of level-set function $\Phi(x)= \|x-(0.5,0.5)^T\|_2-0.25$ on $\hat\Omega=[0,1]^2$}
  \label{fig:levelset}
\end{figure}
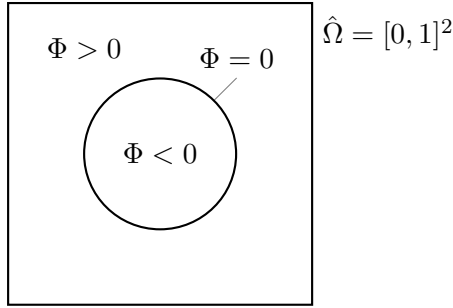

Let $\hat\Omega\subset\RR^d$ denote a polygonal outer domain of dimension $d\in\NN$.
A subdomain $\Omega\subset\hat\Omega\subset\RR^d$ is described using a
level-set function $\Phi$, see \eqref{eq:levelset},
like illustrated in Figure \ref{fig:levelset}.

The boundary $\partial\Omega$ corresponds to the zero level-set $\Phi^{-1}(0) := \{ x \in \hat \Omega \mathrel{\vert} \Phi(x) = 0 \}$.
We call $\partial\Omega$ the interface.
Without loss of generality we restrict the threshold value to 0.
A level-set function with an arbitrary threshold $\alpha\in\RR$ can be obtained by using $\tilde\Phi(x):=\Phi(x)-\alpha$.

Note that it is also possible to describe the motion of a moving interface from an Eulerian point of view in terms of a level-set function and an associated PDE, where $\Phi(x,t)$ satisfies the level-set advection equation
\begin{equation*}
  \label{eq:level_set_advection_equation}
  \Phi_t + \mathbf{v} \cdot \nabla \Phi = 0 \qquad \text{in } \hat\Omega,
\end{equation*}
where $\mathbf{v}(x,t)$ is a velocity field corresponding to the evolution of $\Omega$ and $\Gamma$.

In practice the level-set is given as a discrete scalar function $\Phi_h$.
A common choice is to use a conforming first order discretization,
which is also assumed in this paper.

We consider a mesh $\M(\hat\Omega)$, which is a discretization of the domain $\hat\Omega\subset \RR^d$.
The following properties hold:
$\M(\hat\Omega) = \{\cell_i| i \in \lbrace 0,\dots,N-1\rbrace \}$ forms a partition of $\hat\Omega$ with $\cell_i \cap \cell_j = \emptyset$ if $i \ne j$ and $\bigcup_i \bar \cell_i = \bar{\hat\Omega}$.
$E_i$ denotes the codimension 0 elements of $\M$, i.e. the mesh cells.
In 3D we only consider meshes with cubes, simplices, prisms and pyramids.
In many cases a mesh consists either of cubes or simplices.
In 2D only cubes and simplices are considered.

Apart from mesh cells we identify the following sub entities of $\M$:
faces $\face$ ~(codimension 1), edges $\edge$ ~(codimension $d-1$) and vertices $\vertex$ ~(codimension $d$).
When looking at the codimension it becomes obvious, that in 2D faces and edges are the same.
We still want to distinguish these to avoid confusion, when it comes to generalizations of the algorithm.
For a formal definition of the grid and its properties, as we consider
it, we refer to \cite{Bastian2008}.

The discrete level-set is given as a scalar, continuous and piecewise
multi-linear function.
We choose a representation using vertex values. As $\Phi_h$ is linear
along edges $\edge_i$ it can be written as:
\begin{align*}
  \Phi_h (\vertex_0+t(\vertex_1-\vertex_0))=\Phi_h(\vertex_0)+t\Phi_h(\vertex_1)~,
\end{align*}
with $\vertex_0$ and $\vertex_1$ denoting the vertices associated with $\edge_i$.

\subsection{Numerical integration over implicitly described domains}
\begin{figure}
  \centering
  \parbox{0.25\linewidth}{\centering 2D cube}\quad
  \parbox{0.25\linewidth}{\centering 3D simplex}\quad
  \parbox{0.3\linewidth}{\centering 3D cube}\\
  \parbox{0.25\linewidth}{\includegraphics[height=\linewidth]{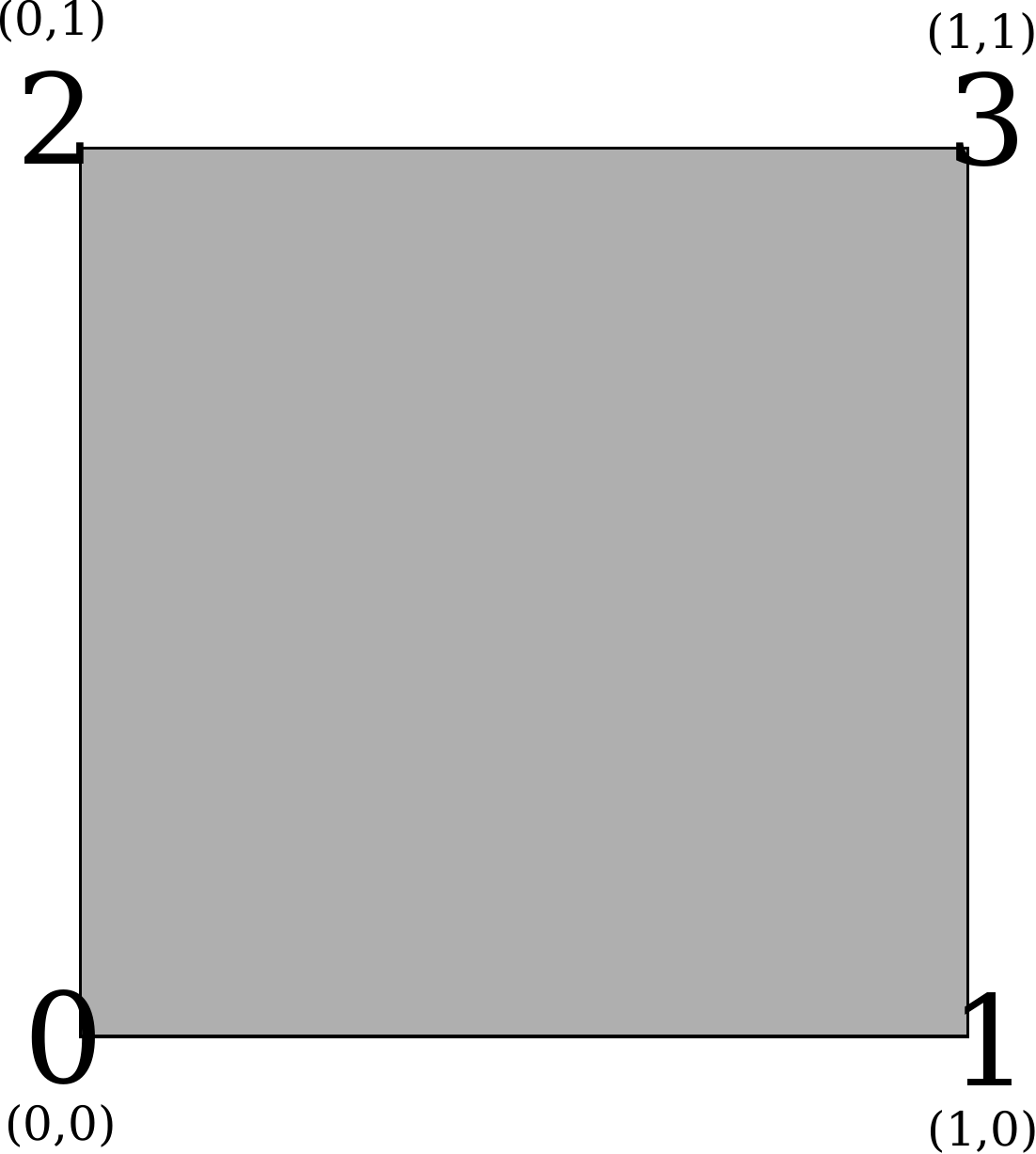}}\quad
  \parbox{0.25\linewidth}{\includegraphics[height=\linewidth]{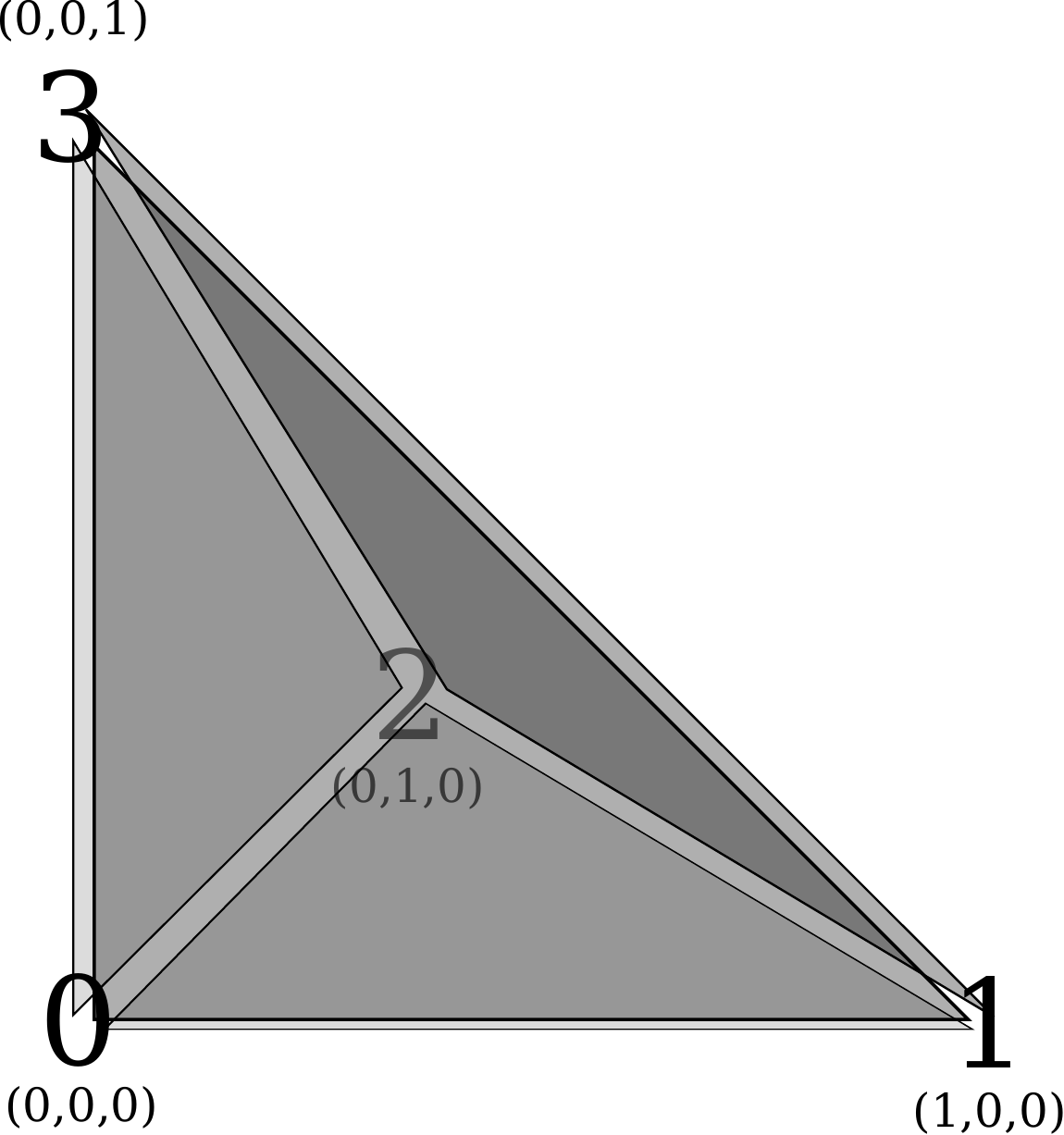}}\quad
  \parbox{0.3\linewidth}{\includegraphics[height=\linewidth]{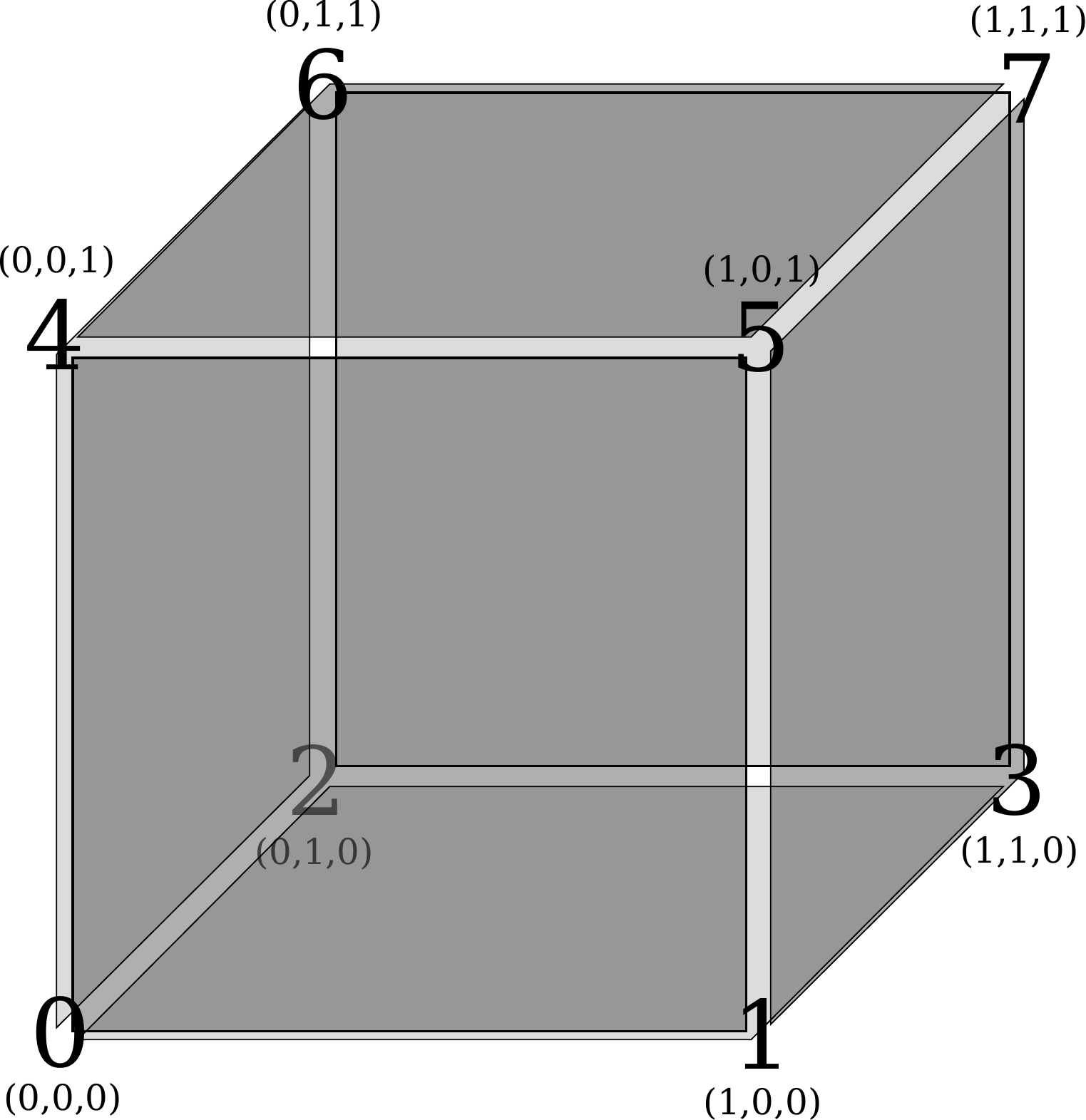}}
  \caption{Vertex numbering and positions in the reference elements.}
  \label{fig:reference_elements}
\end{figure}

The aim of the presented algorithm is to compute the integral of a
function $f:\hat\Omega\to\RR$ over the domain $\Omega$ described by a
discrete level-set function $\Phi_h$. In the following discussion
geometry described by $\Phi_h$ is considered the exact geometry.
The integration over the domain
boundary $\partial\Omega$ or the complement $\hat\Omega\setminus\overline\Omega$ can
be done in a similar manner. Using a triangulation $\M(\hat\Omega)$
of $\hat\Omega$, we compute a sub-triangulation $\{E_i^j\}$ of the intersection of
each element $E_i$ and $\Omega$.
\begin{align*}
  \int_{\lbrace \Phi<0\rbrace}fdx
  & = \sum_{i}\int_{E_i\cap\lbrace \Phi<0\rbrace}fdx \\
  & \approx \sum_{i}\int_{E_i\cap\lbrace \Phi_h<0\rbrace}fdx \\
  & \approx \sum_{i}\sum_{j}\int_{E_i^j}fdx
\end{align*}
The first approximation is unavoidable, the approximation in the last
step is due to our numerical integration and is subject of the
discussion in the section.
The integration over a simple element $E$ (e.g. a simplex or cube) can be done on its reference element using a common quadrature rule:
\begin{align*}
  \int_E fdx = \sum_{k} f(T_E(q_k)) w_k |\det(J_{T_E}(q_k))|
\end{align*}
Here, $T_E$ denotes a bijective map between $E$ and its reference
element and $J_{T_E}$ denotes its Jacobian. In the literature different
choices of reference elements are considered. We follow the definition
used by the DUNE project. A list of these reference elements is given
in Figure \ref{fig:reference_elements}.

\subsection{Topological guarantees}
In order to allow the algorithm to be employed in a wide range of
different methods, we demand that a sub-triangulation fulfills certain topological guarantees.
\begin{enumerate}
\item The connectivity pattern of the cell vertices must be preserved within each subentity.
  In particular this means that vertices connected along an edge, face or volume, should still be connected via the same subentity.
\item This requirement is related to the previous point.
  The curved interface \mbox{$\partial\Omega \equiv \Phi_h^{-1}(0)$} partitions each grid cell into patches belonging to either $\Omega$ or the complement $\hat\Omega\setminus\Omega$.
  We require that the number of patches and their domain association is preserved for the polygonal reconstruction.
\item The vertices of the reconstructed interface lie on the exact
  zero level-set (of the discrete level-set function), i.e. for each vertex $\vertex_i$ the property
  \begin{equation*}
    \Phi_h(\vertex_i) = 0
  \end{equation*}
  holds.
\end{enumerate}

\section{Existing approaches}
\label{sec:existing-approaches}

In this section we discuss existing approaches
for computing a polygonal reconstruction of the
interface $\partial\Omega$ on hexahedral meshes. Extensions to other
element types and other dimensions are always possible.

A well known algorithm is the marching cubes algorithm \cite{Lorensen1987}.
In \cite{engwer2009udg} we described how integration rules for $\Omega$ can be constructed, following the ideas
of the marching cubes algorithm.
In the following we refer to this algorithm as \emph{\MC{}}.
The general idea is to only consider the sign of the function values
in the vertices. Depending on these values $2^{2^\text{dim}}$
different configuration can be distinguished, which can be reduced to
15 basic cases in 3D by taking rotation and mirroring into account. This allows to quickly reconstruct
the interface using predefined rules.

A similar approach was presented in \cite{Min2007}, where the authors basically propose integration rules constructed similar to marching tetrahedrons \cite{gueziec1995mt};
cubes are split into simplices using the Kuhn triangulation \cite{freudenthal1942} and then tetrahedra are triangulated based on their vertex values.
In the following we refer to this algorithm as \emph{\KMT{}}.

As described in different publications, the classic marching cubes algorithm is not able to preserve the topology of the original iso-surface.
\cite{Chernyaev1995} describes how in the marching cubes 33 additional evaluations in the face or volume help to resolve those ambiguities.
A full implementation was published in \cite{Lewiner2003}.
In the following we refer to this algorithm as \emph{\MCxx{}}.
\begin{figure}
  \centering

  \includegraphics[width=0.8\textwidth]{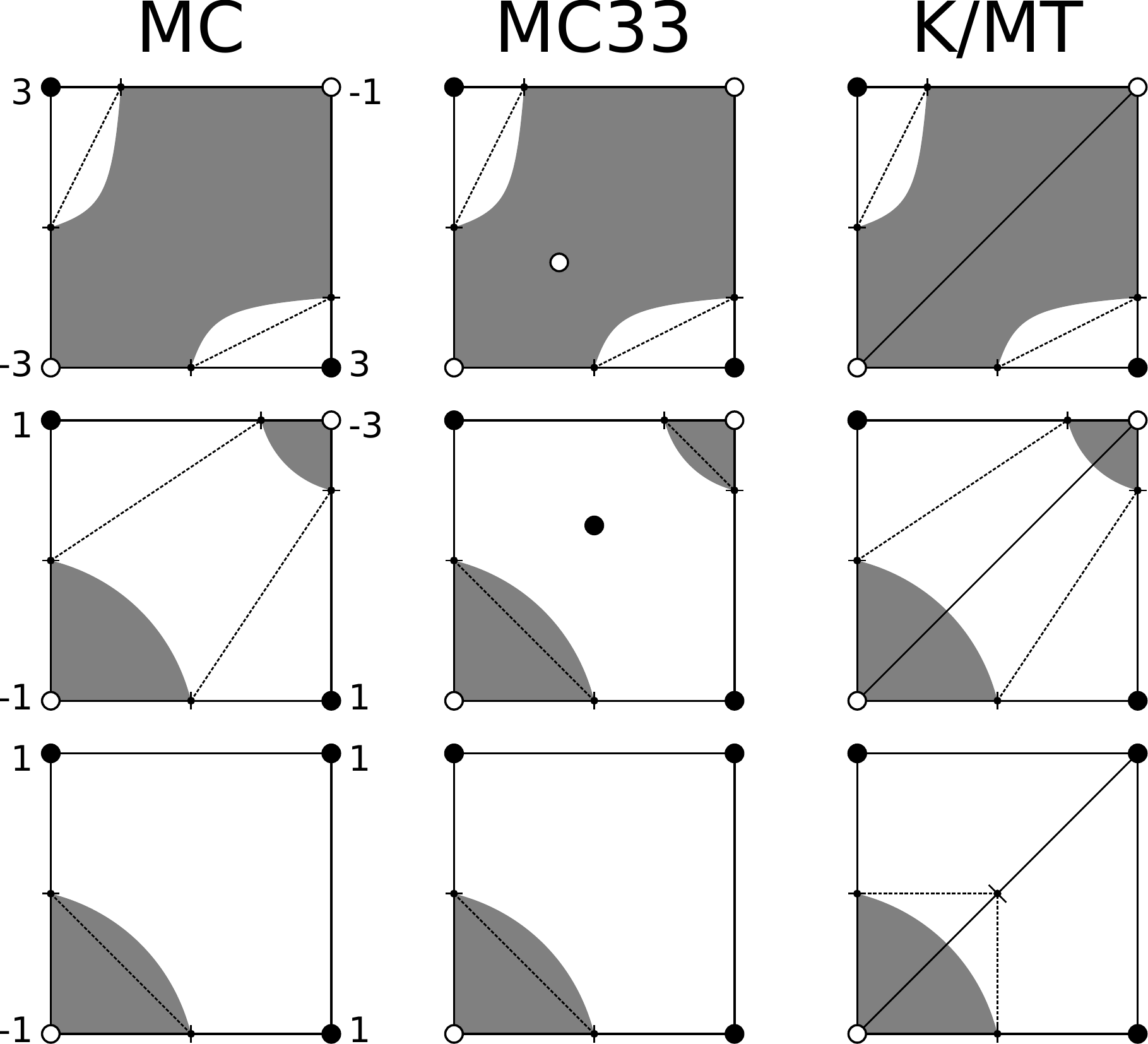}

  \caption{Comparing the existing algorithms: Problems in maintaining
    the topology for the reconstructed interface.
    Gray area: original sub-domain.
    Dots: reconstruction points (black: outside, white: inside).
    Dashed line: reconstructed interface. }
  \label{fig:topologicproblems}
\end{figure}
None of the existing approaches can fulfill all of the afore described
topological guarantees.

While the first property is not fulfilled by the classic \MC{} algorithm, the \MCxx{} algorithm describes an approach, which guarantees this consistency, at least for most cases.
In Section \ref{sec:algorithm} we will describe the cases where the \MCxx{} reconstruction is not sufficient.
The \KMT{} algorithm on the other hand seems to solve the problem of
such ambiguities, as the reconstruction on a simplex is always
unique. The reconstruction is unique, but in general not correct, as
the algorithm implicitly chooses one
solution, due to the previous Kuhn triangulation and by this introduces a preferential direction.

Difficulties arise regarding the third property, which is fulfilled by neither the \MCxx{} algorithm, nor the \KMT{} approach.
In the case of \MCxx{}, nodes in the cube are introduced as weighted
averages of the position of cell vertices to describe a reconstruction
in certain special cases. Such cases can not arise in 2D, but only
in 3D.
Using these nodes in the polygonal reconstruction of the interface
$\partial\Omega$, is sufficient from a visualization point of view,
but not for simulation purposes, as in general the nodes do not lie on
the interface. Depending on the numerical algorithm and the model this
can lead to systematic errors, e.g. in a curvature reconstruction.

\begin{figure}
  \centering
  \begin{tikzpicture}[scale=4]
\fill[lightgray] (.5,.5) circle (0.45);
\draw[gray,step=0.25] (0,0) grid (1,1);
\draw[gray] (0.0,0.0) -- (0.25,0.25);
\draw[gray] (0.0,0.25) -- (0.25,0.5);
\draw[gray] (0.0,0.5) -- (0.25,0.75);
\draw[gray] (0.0,0.75) -- (0.25,1.0);
\draw[gray] (0.25,0.0) -- (0.5,0.25);
\draw[gray] (0.25,0.25) -- (0.5,0.5);
\draw[gray] (0.25,0.5) -- (0.5,0.75);
\draw[gray] (0.25,0.75) -- (0.5,1.0);
\draw[gray] (0.5,0.0) -- (0.75,0.25);
\draw[gray] (0.5,0.25) -- (0.75,0.5);
\draw[gray] (0.5,0.5) -- (0.75,0.75);
\draw[gray] (0.5,0.75) -- (0.75,1.0);
\draw[gray] (0.75,0.0) -- (1.0,0.25);
\draw[gray] (0.75,0.25) -- (1.0,0.5);
\draw[gray] (0.75,0.5) -- (1.0,0.75);
\draw[gray] (0.75,0.75) -- (1.0,1.0);
\draw[thick] (0.25,0.132647574033) -- (0.181801948466,0.181801948466);
\draw[thick] (0.132647574033,0.25) -- (0.181801948466,0.181801948466);
\draw[thick] (0.132647574033,0.25) -- (0.088196601125,0.338196601125);
\draw[thick] (0.05,0.5) -- (0.088196601125,0.338196601125);
\draw[thick] (0.05,0.5) -- (0.0853553390593,0.585355339059);
\draw[thick] (0.132647574033,0.75) -- (0.0853553390593,0.585355339059);
\draw[thick] (0.132647574033,0.75) -- (0.25,0.867352425967);
\draw[thick] (0.5,0.05) -- (0.338196601125,0.088196601125);
\draw[thick] (0.25,0.132647574033) -- (0.338196601125,0.088196601125);
\draw[thick] (0.5,0.95) -- (0.414644660941,0.914644660941);
\draw[thick] (0.25,0.867352425967) -- (0.414644660941,0.914644660941);
\draw[thick] (0.75,0.132647574033) -- (0.585355339059,0.0853553390593);
\draw[thick] (0.5,0.05) -- (0.585355339059,0.0853553390593);
\draw[thick] (0.75,0.867352425967) -- (0.661803398875,0.911803398875);
\draw[thick] (0.5,0.95) -- (0.661803398875,0.911803398875);
\draw[thick] (0.75,0.132647574033) -- (0.867352425967,0.25);
\draw[thick] (0.867352425967,0.25) -- (0.914644660941,0.414644660941);
\draw[thick] (0.95,0.5) -- (0.914644660941,0.414644660941);
\draw[thick] (0.95,0.5) -- (0.911803398875,0.661803398875);
\draw[thick] (0.867352425967,0.75) -- (0.911803398875,0.661803398875);
\draw[thick] (0.867352425967,0.75) -- (0.818198051534,0.818198051534);
\draw[thick] (0.75,0.867352425967) -- (0.818198051534,0.818198051534);
\end{tikzpicture} \quad
  \begin{tikzpicture}[scale=4]
\fill[lightgray] (.5,.5) circle (0.45);
\draw[gray,step=0.25] (0,0) grid (1,1);
\draw[thick] (0.25,0.132647574033) -- (0.132647574033,0.25);
\draw[thick] (0.132647574033,0.25) -- (0.05,0.5);
\draw[thick] (0.05,0.5) -- (0.132647574033,0.75);
\draw[thick] (0.132647574033,0.75) -- (0.25,0.867352425967);
\draw[thick] (0.25,0.132647574033) -- (0.5,0.05);
\draw[thick] (0.25,0.867352425967) -- (0.5,0.95);
\draw[thick] (0.5,0.05) -- (0.75,0.132647574033);
\draw[thick] (0.5,0.95) -- (0.75,0.867352425967);
\draw[thick] (0.75,0.132647574033) -- (0.867352425967,0.25);
\draw[thick] (0.867352425967,0.25) -- (0.95,0.5);
\draw[thick] (0.95,0.5) -- (0.867352425967,0.75);
\draw[thick] (0.867352425967,0.75) -- (0.75,0.867352425967);
\end{tikzpicture} 
  \caption{reconstructed interface of a circle using the \KMT{} approach (left) and \MCxx{} (right)}
  \label{fig:kmt_circle}
\end{figure}
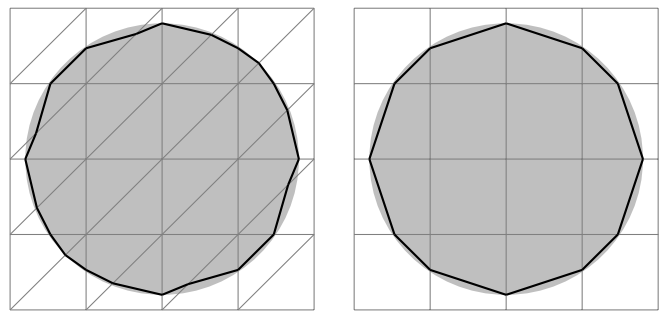

The \KMT{} algorithm first applies a simplex triangulation.
On each simplex, the level-set function is assumed to be linear, which is not correct, as it is actually bi-linear in 2D or tri-linear in 3D.
This leads to reconstructed points, which are slightly off the interface.
In addition, the triangulation into simplices introduces an anisotropy
in the reconstructed interface, see Figure \ref{fig:kmt_circle}.

Further examples of violated consistencies are shown in Section
\ref{sec:algorithm:consistencies}, where we will also discuss how to
prevent these problems using the proposed \ourMC{} algorithm.

\section{Algorithm}
\label{sec:algorithm}

In the following section, we will describe the algorithm for computing
the sub-triangulation of a single grid cell using the vertex values of
the level-set function. All computations can be performed on the
reference element in cell local coordinates.

\subsection{General Idea}

\begin{figure}
  \centering
  \begin{tikzpicture}
  \tikzstyle{inout}=[rectangle, draw, ultra thick]
  \tikzstyle{method}=[rounded corners,draw]
  \tikzstyle{closer}=[node distance=0.8cm]
  \tikzstyle{close}=[node distance=1cm]
  \tikzstyle{far}=[node distance=2cm]
  \tikzstyle{connection}=[->, >=stealth']
  \tikzstyle{entryexit}=[->, >=stealth', thick]
  \node (in) [inout] {\textit{input}};
  \node (entry) [above of=in, closer] {};
  \node (key)[method,close, below of=in] {compute key};
  \node (tests) [method,far, below right of=key] {perform tests};
  \node (coord) [method,far, below left of=tests] {compute coordinates};
  \node (out) [inout,close, below of=coord] {\textit{output}};
  \node (exit) [below of=out, closer] {};
  \path (entry) edge[entryexit] (in)
        (in) edge[connection] (key)
        (key) edge[connection]             node[right] {\textit{ambiguous}} (tests)
              edge[connection, bend right] node[left]  {\textit{non-ambiguous}} (coord)
        (tests) edge[connection] (coord)
        (coord) edge[connection] (out)
        (out) edge[entryexit] (exit);
\end{tikzpicture}
  \caption{structure of the general algorithm for the generation of a sub-triangulation}
  \label{fig:gen-subtriang}
\end{figure}
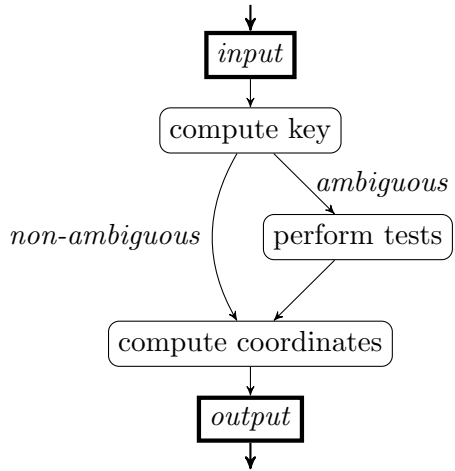

The main idea is to map a set of corner values to a table-index.
The index gives access to different tables which store precomputed general triangulations of the interior, exterior and interface of $\Omega$.
Based on this general triangulation we can compute the actual triangulation by incorporating the corner values.
The basic algorithmic structure is sketched in Figure \ref{fig:gen-subtriang}.

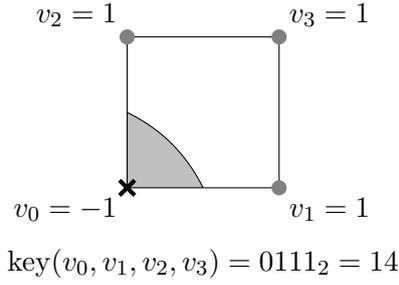
\begin{figure}
  \centering
  \newcommand{\NodeRadius}{0.05}
\begin{tikzpicture}[scale=2]
  \begin{scope}[domain=0:.5, samples=21]
    \fill[lightgray] (0,0) -- plot (\x,{1.0+1.0/(2.0*\x-2.0)});
    \draw (0,0) -- plot (\x,{1.0+1.0/(2.0*\x-2.0)});
  \end{scope}
  \draw (0,0) rectangle (1,1);
  \foreach \x/\y/\l/\a in {1/0/$v_1=1$/below right, 0/1/$v_2=1$/above left, 1/1/$v_3=1$/above right} {
    \fill[gray] (\x,\y) circle (\NodeRadius) node[\a, black] {\l};
  }
  \foreach \x/\y/\l/\a in {0/0/$v_0=-1$/below left} {
    \draw[ultra thick] (\x,\y) -- +(\NodeRadius,\NodeRadius);
    \draw[ultra thick] (\x,\y) -- +(-\NodeRadius,\NodeRadius);
    \draw[ultra thick] (\x,\y) -- +(\NodeRadius,-\NodeRadius);
    \draw[ultra thick] (\x,\y) -- +(-\NodeRadius,-\NodeRadius);
    \draw (\x,\y) node[\a] {\l};
  }
  \node (anker) at (0.5,0) {};
  \node[below of=anker] {$\operatorname{key}(v_0,v_1,v_2,v_3)=0111_2=14$};
\end{tikzpicture}
  \caption{generation of the key for a two dimensional square domain based on
           the level set values in the corners}
  \label{fig:key-generation}
\end{figure}

The generation of the table index is subdivided into two steps: the computation of a key and, if necessary, the application of certain tests if the key is ambiguous.
The key of a sequence of corner values $(v_0,\dots,v_{n-1})$ of
$\Phi_h$ is defined as
\begin{align*}
  \operatorname{key}(v_0,\dots,v_{n-1}):=\sum_{i=0}^{n-1}\chi_{\RR^{\ge0}}(v_i)2^i,
\end{align*}
where $\chi_M$ denotes the indicator function
\begin{align*}
  \chi_M(x):=
  \begin{cases}
    1 &,x\in M \\
    0 &,x\not\in M
  \end{cases}.
\end{align*}
An example for the key computation on a two dimensional square is shown in Figure \ref{fig:key-generation}.
For some cases, this key describes the topology completely.
Other cases require additional tests to determine the triangulation.
We distinguish between two types of tests: face tests and center tests.

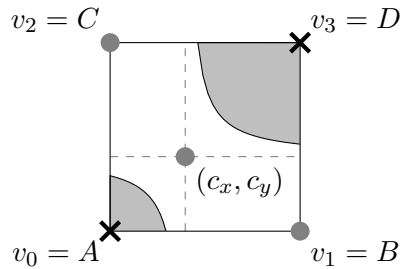
\begin{figure}[b]
  \centering
  \newcommand{\NodeRadius}{0.05}
\newcommand{\valueA}{0.5}
\newcommand{\valueBC}{1.2}
\newcommand{\valueD}{1.4}
\begin{tikzpicture}[scale=2.5]
  \begin{scope}[domain=0.0:{\valueA/(\valueBC+\valueA)}, samples=21]
    \fill[lightgray] (0,0) -- plot (\x,{(\valueA-(\valueBC+\valueA)*\x)/(\valueBC+\valueA-(\valueA+2.0*\valueBC+\valueD)*\x)});
    \draw (0,0) -- plot (\x,{(\valueA-(\valueBC+\valueA)*\x)/(\valueBC+\valueA-(\valueA+2.0*\valueBC+\valueD)*\x)});
  \end{scope}
  \begin{scope}[domain={\valueBC/(\valueBC+\valueD)}:1.0, samples=21]
    \fill[lightgray] (1,1) -- plot (\x,{(\valueA-(\valueBC+\valueA)*\x)/(\valueBC+\valueA-(\valueA+2.0*\valueBC+\valueD)*\x)});
    \draw (1,1) -- plot (\x,{(\valueA-(\valueBC+\valueA)*\x)/(\valueBC+\valueA-(\valueA+2.0*\valueBC+\valueD)*\x)});
  \end{scope}
  \draw[dashed, gray] ({(\valueA+\valueBC)/(\valueA+2.0*\valueBC+\valueD)},0) -- ({(\valueA+\valueBC)/(\valueA+2.0*\valueBC+\valueD)},1.0);
  \draw[dashed, gray] (0,{(\valueA+\valueBC)/(\valueA+2.0*\valueBC+\valueD)}) -- (1.0,{(\valueA+\valueBC)/(\valueA+2.0*\valueBC+\valueD)});
  \fill[gray] ({(\valueA+\valueBC)/(\valueA+2.0*\valueBC+\valueD)},{(\valueA+\valueBC)/(\valueA+2.0*\valueBC+\valueD)}) circle (\NodeRadius) node[below right,black] {$(c_x,c_y)$};
  \draw (0,0) rectangle (1,1);
  \foreach \x/\y/\l/\a in {1/0/$v_1=B$/below right, 0/1/$v_2=C$/above left} {
    \fill[gray] (\x,\y) circle (\NodeRadius) node[\a, black] {\l};
  }
  \foreach \x/\y/\l/\a in {0/0/$v_0=A$/below left, 1/1/$v_3=D$/above right} {
    \draw[ultra thick] (\x,\y) -- +(\NodeRadius,\NodeRadius);
    \draw[ultra thick] (\x,\y) -- +(-\NodeRadius,\NodeRadius);
    \draw[ultra thick] (\x,\y) -- +(\NodeRadius,-\NodeRadius);
    \draw[ultra thick] (\x,\y) -- +(-\NodeRadius,-\NodeRadius);
    \draw (\x,\y) node[\a] {\l};
  }
\end{tikzpicture}
  \caption{Interface of an ambiguous case on a quadratic face. $v_0$ and $v_3$ are negative and $v_1$ and $v_2$ are positive. The interface forms a hyperbola with its center located at $(c_x, c_y)$.}
  \label{fig:hyperbola}
\end{figure}

The triangulation of a quadratic face is ambiguous, if diagonally opposing vertices have the same and neighboring vertices a different sign.
For such a case, a test needs to be performed to decide if the face center is in $\Omega$ or $\hat\Omega\setminus\overline\Omega$.
The interface on the ambiguous face builds a hyperbola.
An example for such an interface can be seen in Figure \ref{fig:hyperbola}.
The sign of the face center is given as the sign of the hyperbolas center.
For vertex values $A,B,C,D$ at vertices $0,1,2,3$
respectively, the center is located (in local
coordinates of the face) at
 \begin{align*}
  (c_x,c_y):=\left(\frac{A-C}{A-B-C+D},
    \frac{A-B}{A-B-C+D}\right).
\end{align*}
and its sign can be computed as
\begin{align*}
  \sgn(A)\sgn(AD-BC).
\end{align*}
For further details, we refer to \cite{Chernyaev1995}. In three space
dimensions further ambiguities can occur.

Vertices which are not connected via a face can still be connected
through the volume. This means that the iso-surface forms a kind of
tube through the element.
A method to test this has been presented in \cite{Chernyaev1995} which we will recall here.
Let $v$ and $w$ denote two vertices on opposite ends of a diagonal of the cube
which have the same sign and are not connected via a set of faces. If and only
if $v$ and $w$ are connected through the cubes volume, there exists a plane
parallel to a face on which the projections of $v$ and $w$ are connected. This
connection can be explicitly calculated in terms of the cubes vertex values.
Let $A_0,B_0,C_0$ and $D_0$ denote the values at the vertices of the face containing $v$
and $A_1,B_1,C_1$ and $D_1$ denote the values at the vertices of the face containing $w$,
while $A_0$ denotes the value at $v$ and $D_1$ denotes the value at $w$.
We assume that $A_0$ and $D_1$ have a positive sign.
The value along an edge perpendicular to the two faces is given by the linear functions $A_t,B_t,C_t,D_t$ for $t\in[0,1]$.
To calculate the connection, we consider the following quadratic function
\begin{equation}
  p:[0,1]\to\RR; t\mapsto A_tD_t-B_tC_t=at^2+bt+c\,, \label{eq:quadratic}
\end{equation}
where the coefficients $a,b,c\in\RR$ are given by
\begin{align*}
  a =& (A_1-A_0)(D_1-D_0)-(C_1-C_0)(B_1-B_0) \\
  b =& D_0(A_1-A_0)+A_0(D_1-D_0)-B_0(C_1-C_0)-C_0(B_1-B_0) \\
  c =& A_0D_0-B_0C_0
\end{align*}
The value of $p$ corresponds to the value of the center of the hyperbola on the given plane.
The positive areas are connected, if $p$ has a maximum at $t_{\operatorname{max}} \in [0,1]$ with $p(t_{\operatorname{max}})>0$ and the values $A_t,B_t,C_t,D_t$ have the correct sign. Otherwise they are separated.
Again, we refer to \cite{Chernyaev1995} for a detailed description.

The necessary tests for a specific case form a tree structure, which
is traversed in order to calculate the table index.
A table containing the necessary tests is available in an online resource of \cite{Lewiner2003}.
Note though, that we found two tests to be missing: for cases 10 and 12 with positive tests for both ambiguous faces, one has to check whether the exterior vertices are connected through the cube.

Once the table index has been computed, it can be used to retrieve general triangulations.
Such a triangulation describes its local elements using the vertex numbering of the reference element.
The actual coordinates are computed using the vertex values.
We distinguish between five types of vertices which are described in Section \ref{sec:construction_of_vertices}.

\subsection{Additional consistencies}
\label{sec:algorithm:consistencies}

As noted in Section \ref{sec:preliminaries} additional consistencies
are required to ensure the topological guarantees. The \MCxx{}
improves the classical \MC{} algorithm to ensure topological
consistency for the reconstructed interface. While this is sufficient
for visualization, it is in general not sufficient for numerical
simulations.

\begin{figure}[t]
 \centering
 \includegraphics[width=0.3\textwidth]{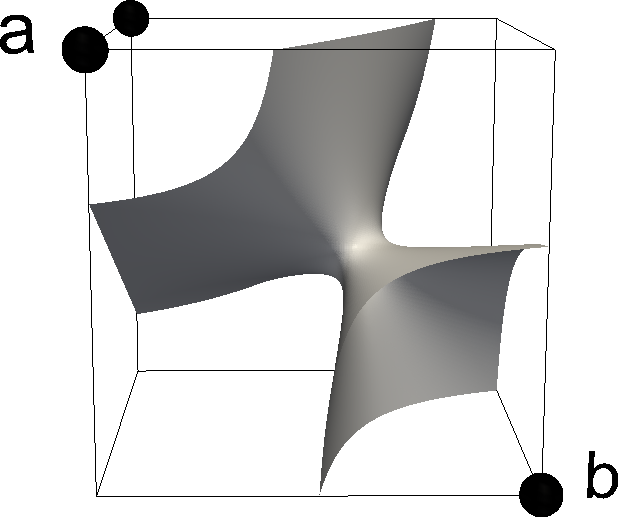}
 \includegraphics[width=0.3\textwidth]{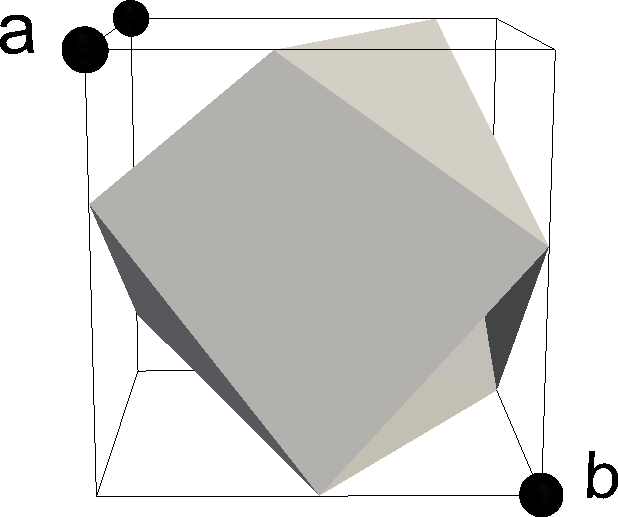}
 \includegraphics[width=0.3\textwidth]{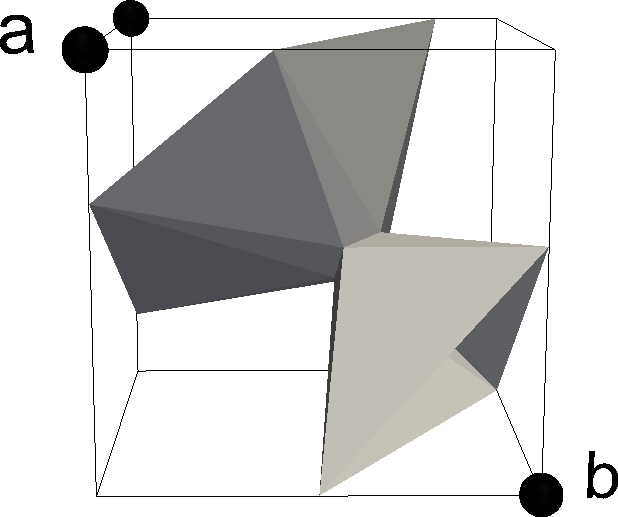}
 \caption{
    Violation of the topologic consistency in a face (front), using \MCxx{} and the correct reconstruction using \ourMC.
    Left: discrete level-set, center: reconstruction from \MCxx{}
    (case 6.1.2 of \cite{Lewiner2003}), right: \ourMC{}
    reconstruction.
    For the choice of vertex values see Table \ref{tab:face_error}.}
   \label{fig:face_error}
\end{figure}

\begin{table}[t]
  \centering
    \begin{tabular}{l|rrrrrrrr}
      \hline\hline
      \rule{0pt}{3ex}Vtx & 0 & 1 & 2 & 3 & 4 & 5 & 6 & 7 \\
      Pos & (0,0,0) & (0,0,1) & (0,1,0) & (0,1,1) & (1,0,0) & (1,0,1) & (1,1,0) & (1,1,1) \\
      Value & -4 & 4 & -1 & -1 & 2 & -3 & 5 & -1\\[1ex]
      \hline\hline
    \end{tabular}
  \caption{Vertex values for the example illustrated
    in Figure \ref{fig:face_error}.}\label{tab:face_error}
\end{table}

We want to illustrate the issue with a concrete example. Figure
\ref{fig:face_error} shows on the left the correct interface of a
tri-linear function with vertex values as listed in table
\ref{tab:face_error}. After resolving the ambiguities the \MCxx{}
algorithm provides a reconstruction where parts of the interface are
not in the interior of the cell anymore, but are pushed into the
frontal face of the cube; the volume which lies in-between vanishes when using
the discrete level-set function. This leads to an
unresolvable inconsistency between the trace of the volume
reconstruction and a direct reconstruction in the cube face.
Therefore, these reconstructions do not fulfill the topological
guarantees and are not accurate enough for numerical purposes.  This
issue can be resolved, by adding interior points to reconstruct the
interface and especially the vanishing volume more accurately.

\subsection{Construction of Vertices}
\label{sec:construction_of_vertices}

\begin{figure}[t]
 \centering
 \includegraphics[width=0.3\textwidth]{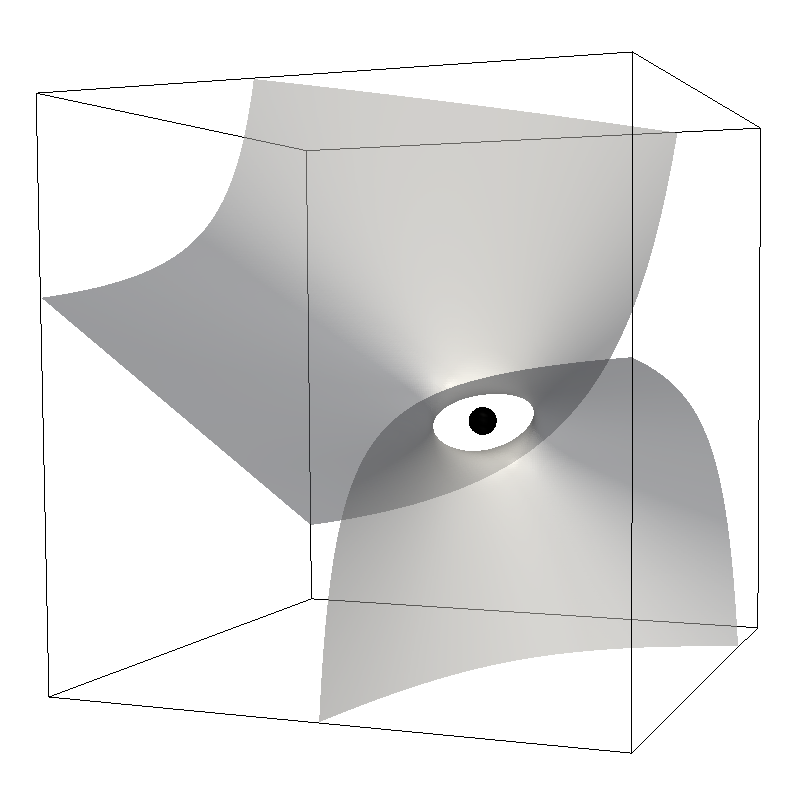}
 \includegraphics[width=0.3\textwidth]{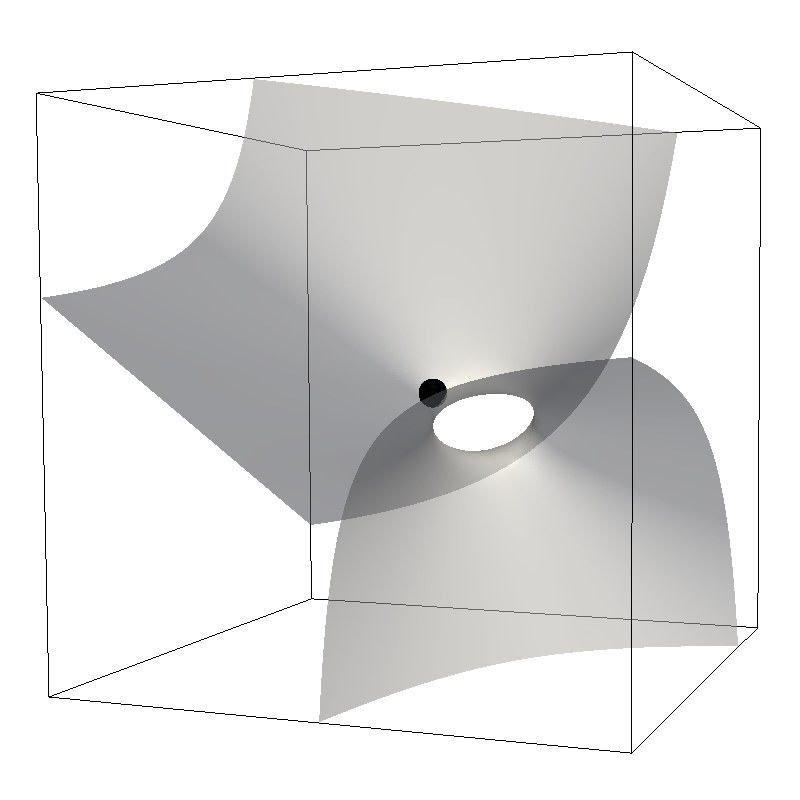}
 \includegraphics[width=0.3\textwidth]{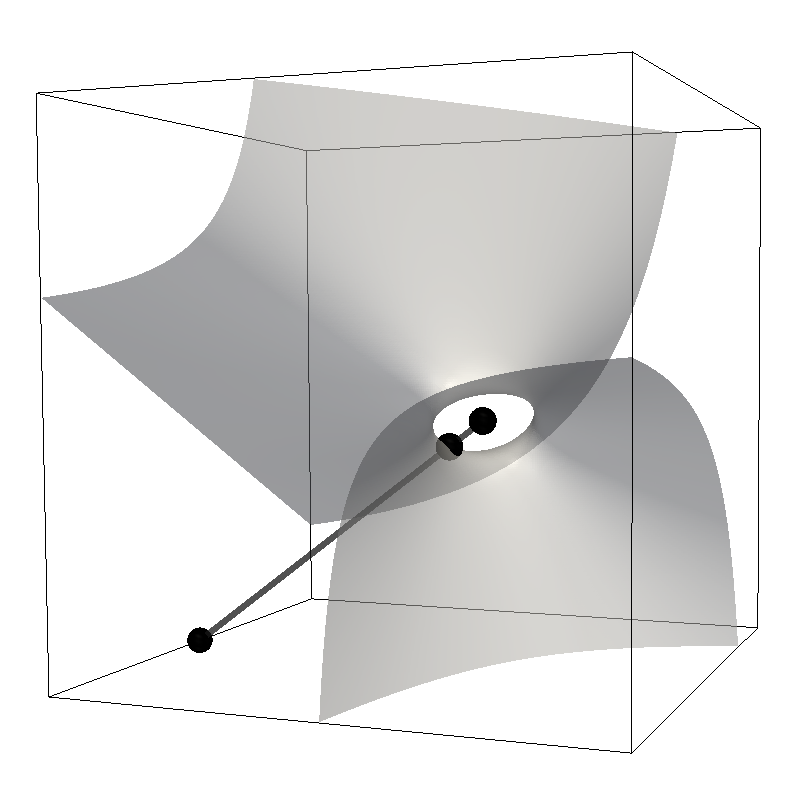}
 \caption{
   Maximum vertex (left), root vertex (center) and intersection vertex (right) for the basic case 6.1.2. For the choice of vertex values see Table \ref{tab:face_error}.  }
   \label{fig:vertex_max_and_root}
\end{figure}

Based on the general reconstruction rule stored in the look-up table
the actual reconstruction is computed. This requires the construction
of vertices based on the actual vertex values. We distinguish between
five types of vertices:
\begin{enumerate}
\item \emph{original vertices}:\\
  For the volume reconstruction we have to incorporate the vertices of the
  reference element, our implementation is based on a dimension
  independent numbering used in the DUNE framework (see Figure \ref{fig:reference_elements}).
\item \emph{edge vertices}:\\
  These vertices describe the position of the interface along an edge
  of the reference element. As the level-set is linear along an edge
  the position is easily computed as weighted average
  \begin{equation*}
    x_{\{A,B\}} = x_A \frac{B}{B - A} + x_B \frac{A}{A - B}
  \end{equation*}
  where $A,B$ denote the level-set values of the vertices at $x_A,x_B$ of the edge.
\end{enumerate}
In order to resolve the issues described in Section \ref{sec:algorithm:consistencies} we require additional vertices in the interior of the cell.
\begin{enumerate}[resume]
\item \emph{maximum vertices}:\\
  By maximum vertices we denote the vertices given by the maximum of the quadratic function $p$ defined in equation \eqref{eq:quadratic} of the interior connectivity test.
  Such a vertex is located at the ``center'' of the constriction of a tube-like structure. For an example
  of such a structure see Figure \ref{fig:face_error} and an example of a maximum vertex see Figure \ref{fig:vertex_max_and_root}.
  They are not used directly in the triangulation since they do not lie on the interface, but are used as helper vertices for vertices of type \ref{it:interiorvertices}.
\item \emph{root vertices}: \\
  Root vertices are defined similar to maximum vertices and are given by the roots of the quadratic function $p$ defined in equation \eqref{eq:quadratic}.
  If the interior connectivity test is positive,
  there is a plane (defined by $t_{\operatorname{max}}$) parallel to a given face of the cube where the restricted interface forms a hyperbola.
  Moving the plane along its normal and following the path of the center of the hyperbola, the root vertices are given by the intersection of the path with the interface.
  An example for a root vertex is shown in Figure \ref{fig:vertex_max_and_root}.
  Note that such a point always exists as in this case the sign at the
  center of the hyperbola is known to change.
\item\label{it:interiorvertices} \emph{interior vertices}:\\
  These points are
  computed as the intersection between the iso-surface and an a-priori
  chosen line. As the line is not axis-aligned and the level-set
  function $\Phi_h$ is tri-linear we can have up to three intersection
  points along the line. We compute all intersection simultaneously
  using the Aberth algorithm \cite{aberth1973}.
  The latter is a Newton-type iterative method for finding the roots of a
  polynomial. It makes sure that no two approximations converge to the same
  roots.
\end{enumerate}

\begin{example}[Construction of interior; case 6.1.2]
  \label{example:construction}
  To illustrate the different steps of the algorithm,
  we consider the example illustrated in Figure \ref{fig:face_error}.
  The according vertex values are listed in Table
  \ref{tab:face_error} and the vertex value yield the key $01001010_2=82$,
  which correspond to the basic case 6 (see \cite{Lewiner2003}). This
  case is ambiguous as vertices $1$ and $6$, which are connected via a diagonal through the cube,
  have the same sign and are not connected over faces of the cube.
  To perform the corresponding test, we compute the coefficients of the quadratic function with respect to the front face, which is given as
  \begin{align*}
    p(t)=-42t^2+37t-8
  \end{align*}
  This function has a maximum at $t_{\operatorname{max}}=7/84\in[0,1]$ and $p(t_{\operatorname{max}})=5/168>0$.
  The values of the projections are $A_{t_{\operatorname{max}}}=1/12$, $B_{t_{\operatorname{max}}}=-1$, $C_{t_{\operatorname{max}}}=-9/14$ and $D_{t_{\operatorname{max}}}=3/14$.
  These values show, that in the plane of $t_{\operatorname{max}}$, $A_0$ and $D_1$ are connected.
  This result corresponds to the basic case 6.1.2.

  For this basic case we query the look-up tables and obtain
  reconstruction rules for interior, exterior and interface, which
  describes a set of parts as a list of vertices and a list of
  elements/faces. These describe primitive geometric objects, suitable
  for numerical integration.
  In our implementation the reconstruction of the interior,
  exterior and interface contain lists of 23, 21 and 15 parts
  respectively.  The vertices in these lists are encoded via their
  associated subentity in the reference element.  The lists for the specific case
  6.1.2 contain for example reference vertex $3$, which is located at
  $(1,1,0)$.  Another example is the intersection vertex between the
  reference vertices $1$ (located at $(1,0,0)$) and $3$.  This vertex
  is found by linear interpolation from the vertex values $4$ and $-1$
  which gives $(1,4/5,0)$.  A more complex example is the intersection
  vertex between the maximum vertex and the middle point on the edge
  of vertices $0$ and $2$ as shown in the right Figure
  \ref{fig:vertex_max_and_root}.  First, we compute the maximum vertex
  by computing the center of the hyperbola in the plane parallel to
  the front face at $t_{\operatorname{max}}=7/84$. This vertex is
  given as approximately $(0.61017, 0.45238, 0.44552)$ (see Figure
  \ref{fig:vertex_max_and_root} left). In general $\Phi$ is locally a
  trilinear function, thus we use the Aberth-method to
  get the intersection with interface of the line between this maximum
  vertex and the center $(0,0.5,0)$ of the edge between vertices $0$
  and $2$. The intersection vertex is approximately $(0.54437,
  0.45752, 0.39748)$.  We compute the vertex position for each vertex
  in the list by the methods defined above.  From these vertices we
  can construct the volume of the interior and exterior as the union
  of the respective parts.
\end{example}

\section{Implementation}
\label{sec:implementation}

An implementation of the algorithm is available via
\texttt{python-pip}, all tests performed in this paper are available
for download from ZENODO~\cite{download-tpmc}. A python
library contains a list of all basic cases and their associated
information. From the basic cases the python code computes the full
list and can then generate C++ code. This generated code only uses
standard C++ components and can easily be used in other projects. The
code is in productive use in DUNE for our implementation of the Unfitted Dicontinuous
Galerkin method~\cite{engwer2010:duneudg}.

The C++ library stores all information in a set of different look-up
tables. The tests are stored as binary trees in breadth-first
order as an implicit data structure in arrays. The library provides a basic interface to perform the
tests and export geometric reconstructions as sets
of sets of vertices. The library is templated, so that it can work
with user defined container types and coordinate types. This allows to incorporate it
into existing projects, without having to change internal data
structures. As the look-up tables are stored as plain C arrays, they can also
be incorporated into C or Fortran applications.

\subsection{Look-up tables}
The geometric information of the sub-triangulation parts for each case is stored in three different tables:
\begin{enumerate}
\item Reconstruction of the interface $\Gamma$
\item Reconstruction of $\Omega$
\item Reconstruction of the complement of $\Omega$
\end{enumerate}
Note that we store two different tables for $\Omega$ and its complement.
Naively a sub-triangulation of the complement
$\hat\Omega\setminus\Omega$ could be computed by simply inverting the sign of the level set function.
However, the resulting triangulation of $\hat\Omega$ will in general
no longer form a partition, it might contain holes or overlapping parts.
To avoid this, both sub-triangulations have to be constructed with a matching interface.
The vertices of the parts are encoded using their position in the reference element.
The ordering of the vertices of each interface part is done such that
the normal vectors to $\Gamma$ can be computed using the right hand
rule, i.e. the direction of the normal vector of an interface part with vertices $x_0,x_1,x_2,\dots$ is given by the cross product of the tangential axes: $(x_1-x_0)\times(x_2-x_0)$.

In addition to the general geometric information we store information about the connectivity pattern w.r.t. each sub-domain.
For each sub-triangulation part and each vertex we store an index of the local connected component it belongs to.
Using this information, the global connectivity pattern can be computed without geometric coordinate comparisons.

Look-up tables are provided for all primitive geometries up to
dimension three, i.e. for lines, triangles, quadrilateral, tetrahedra,
prisms, pyramids and hexahedra. Thus the library can be used with a
broad range of discrete level-set representations.

\subsection{Example}

Following the sketch in Figure \ref{fig:gen-subtriang} we first compute a
key based on given vertex values. If we encounter an ambiguous case,
this involves a set of different tests, which are retrieved from a
binary tree stored in
table \texttt{table\_cube3d\_mc33\_face\_test\_order}. In the next
step we compute all internal vertices listed in; as these computations might be
expensive, it is done in a seperate step to avoid repetitive
computation. The list of necessary
vertices is stored in \texttt{table\_cube3d\_vertices}. As the last step, we retrieve the
desired information, e.g. the reconstruction of the interior domain
$\Omega$. In this particular example
\texttt{table\_cube3d\_cases\_offsets} yields offsets into a
reconstruction table, e.g. \texttt{table\_cube3dsym\_codim\_0\_interior}
for the interior domain.

The corresponding C++ user code to compute the construction described in
Example \ref{example:construction} is listed in Algorithm \ref{alg:usercode}.
\begin{algorithm}[t]
\caption{C++ user code to compute an interior domain decomposition, as
  described in Example \ref{example:construction}.}
\begin{lstlisting}
// some necessary typedefs
static unsigned int dim = 3;
typedef std::array<double,dim> Coordinate;
// instantiate the algorithm and some additional
// structures
tpmc::MarchingCubes<double, dim, Coordinate> mc;
tpmc::GeometryType geometry = tpmc::cube;
std::vector<Coordinate> vertices;
vertices.reserve(mc.getMaximalVertexCount(geometry));
// describe the geometry
std::vector<double> vertex_values =
    {-4.0, 4.0, -1.0, -1.0, 2.0, -3.0, 5.0, -1.0};
// compute the key
std::size_t key = mc.getKey(vertex_values.begin(),
                            vertex_values.end());
// compute additional vertices and append them
mc.getVertices(vertex_values.begin(), vertex_values.end(),
               key,
               std::back_inserter(vertices)
              );
// compute interior reconstruction and store them in
// decomposition
std::vector< std::vector<int> > decomposition;
mc.getElements(geometry, key, tpmc::InteriorDomain,
               std::back_inserter(decomposition)
              );
\end{lstlisting}
\label{alg:usercode}
\end{algorithm}

The decomposition of the interior domain is described as a cell-to-vertex
map in the form $C_i \rightarrow \{ v_j \}$. Each entry describes a
primitive geometry.
This map is stored in the \lstinline!decomposition! variable.
Given the vertex numbers of a cell the actual coordinates are
obtained by a lookup in the vector of vertex coordinates, stored in
\lstinline!vertices!. This description follows the approach used in
most mesh file formats.

\section{Evaluation}
\label{sec:evaluation}

We present different numerical tests to evaluate the accuracy, the
robustness and the efficiency of the proposed algorithm.

\subsection{Surface Integral over two entangled Tori}
To evaluate the convergence rate of the integration method, we use a test level-set function representing two tori.
We set $\hat\Omega=[0,1]^3\subset\RR^3$ as the outer domain.
For a unit vector $n\in\RR^3$ we define its associated hyperplane as $H(n):=\lbrace x\in\hat\Omega|\langle x,n\rangle=0\rbrace$.
The orthogonal projection onto $H(n)$ is then given by $P(n)(x)=x-\langle x,n\rangle n$.
The level-set function of a torus with radius $R\in\RR$ and tube
radius $r\in\RR$ in a hyperplane $H(n)$ is given by
\begin{align*}
 \Phi_{R,r}(n)(x):=\sqrt{(\|P(n)(x)\|-R)^2+\langle x,n\rangle^2}-r~.
\end{align*}%
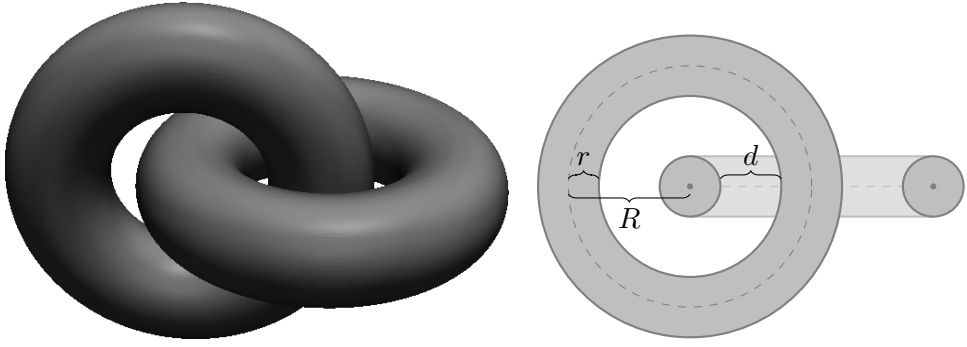
\begin{figure}
 \centering
 \includegraphics[width=0.45\linewidth]{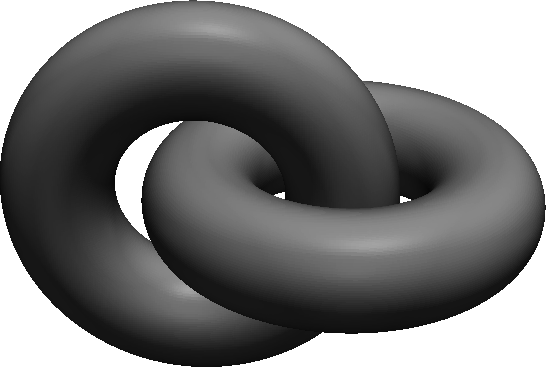}\hspace{0.25cm}
 \begin{tikzpicture}[scale=2]
  \fill[lightgray!50!white] (0,-0.2) rectangle (1.6,0.2);
  \draw[gray!50!white,dashed] (0,0) -- (1.6,0);
  \draw[gray!50!white,thick] (0,0.2) -- (1.6,0.2);
  \draw[gray!50!white,thick] (0,-0.2) -- (1.6,-0.2);
  \fill[lightgray] (0,0) circle (1.0);
  \fill[white] (0,0) circle (0.6);
  \begin{scope}
    \clip (0,0) circle (0.6);
    \fill[lightgray!50!white] (0,-0.2) rectangle (1.6,0.2);
    \draw[gray!50!white,dashed] (0,0) -- (1.6,0);
    \draw[gray!50!white,thick] (0,0.2) -- (1.6,0.2);
    \draw[gray!50!white,thick] (0,-0.2) -- (1.6,-0.2);
  \end{scope}
  \fill[lightgray] (0,0) circle (0.2);
  \fill[lightgray] (1.6,0) circle (0.2);
  \fill[gray] (0,0) circle (0.02);
  \fill[gray] (1.6,0) circle (0.02);
  \draw[gray,dashed] (0,0) circle (0.8);
  \draw[gray,thick] (0,0) circle (1.0);
  \draw[gray,thick] (0,0) circle (0.6);
  \draw[gray,thick] (0,0) circle (0.2);
  \draw[gray,thick] (1.6,0) circle (0.2);
  \draw[decorate,decoration={brace,amplitude=3pt}] 
      (-0.8,0.05) -- (-0.6,0.05); 
  \draw (-0.7,0) node[above=0.15] {$r$};
  \draw[decorate,decoration={brace,amplitude=3pt,mirror}] 
      (-0.8,-0.05) -- (0,-0.05); 
  \draw (-0.4,0) node[below=0.15] {$R$};
  \draw[decorate,decoration={brace,amplitude=3pt}] 
      (0.2,0.05) -- (0.6,0.05); 
  \draw (0.4,0) node[above=0.15] {$d$};
\end{tikzpicture}
 \caption{reconstructed interface of the test problem of two entangled tori and a cut along the hyperplane associated with one torus.}
 \label{fig:tori_interface}
\end{figure}%
For $\phi\in[0,\pi)$, we set
\begin{align*}
 n_1(\phi):=\begin{pmatrix}
       \cos(\phi)\\
       \sin(\phi)\\
       0
      \end{pmatrix},\quad
 n_2(\phi):=\begin{pmatrix}
       -\sin(\phi)\\
       \cos(\phi)\\
       0
      \end{pmatrix},\quad
 n_3:=\begin{pmatrix}
       0\\
       0\\
       1
      \end{pmatrix},\quad
 c:=\begin{pmatrix}
     0.5\\
     0.5\\
     0.5
    \end{pmatrix}
\end{align*}
and define the level-set function of two entangled tori (see Figure
\ref{fig:tori_interface}), centered at $c$:
\begin{align*}
 \Phi_{R,r,\phi}(x):=\min\Big(&\Phi_{R,r}(n_1(\phi))\Big(x-c+\frac{R}{2}\cdot n_3\Big),\\
      &\Phi_{R,r}(n_2(\phi))\Big(x-c-\frac{R}{2}\cdot n_3\Big)\Big)~.
\end{align*}
The surface area of the level-set $\Phi_{R,r,\phi}$ can be computed
analytically as
\begin{align*}
 S(\Phi_{R,r,\phi}):=\int_{\Gamma_{R,r,\phi}}1ds=8\pi^2Rr
\end{align*}
The surface of the reconstructed interface using a structured grid with diameter $h$ is denoted by $S_h(\Phi_{R,r,\phi})$.
The relative error of a reconstructed surface is defined as
\begin{align*}
 e_S:=\frac{|S(\Phi_{R,r,\phi})-S_h(\Phi_{R,r,\phi})|}{S(\Phi_{R,r,\phi})}
\end{align*}

We define $d = R - 2r$ as the distance between the two tori. For $d =
0$ the two tori touch and the topology changes, the aspect ratio $R/r$ in this case is
$2$.

We set $R=0.25$ with a distance $d=0.1$,
i.e. we choose the minor radius $r=\frac{R-d}{2}=0.075$ and set $\phi=0$.
In order to investigate the robustness of the method to grid perturbation, we shift the center of the tori
randomly by a maximum $h$. For each grid size, we generate $30$ random datasets.
The results of the relative surface error for different grid sizes are
listed in Table \ref{tab:plot_accuracy}.

\begin{table}
  \centering
    \begin{tabular}{@{\quad}r@{\qquad\quad}l@{\quad$-$\quad}l@{\quad\qquad}l@{~~$\pm$~~}l@{\quad}}
      \hline\hline
      \rule{0pt}{1.9ex}grid size & $\min(e_S)$ & $\max(e_S)$ & $\overline{e_S}$ & $\sigma_1$ \\\hline
      16  & 0.0451   & 0.0509   & 0.046   & $3.2\%$ \\ 
      32  & 0.0108   & 0.0114   & 0.0110   & $1.8\%$ \\ 
      64  & 0.00266  & 0.00275  & 0.00270  & $1.1\%$ \\ 
      128 & 0.0006672   & 0.00067563 & 0.0006738 & $0.3\%$ \\ 
      256 & 0.0001674 & 0.0001686 & 0.0001681 & $0.2\%$ \\ 
      \hline\hline
    \end{tabular}
  \caption{Relative error $e_S$ of the surface of the two tori test problem computed with \ourMC{} as compared to the analytic solution. The level-set is discretized on a three dimensional structured grid with rectangular element. For each grid size we generate $30$ randomly shifted datasets}
  \label{tab:plot_accuracy}
\end{table}

We can observe a second order convergence rate, as it is expected for a polygonal approximation.
In addition, the computation is robust with respect to grid perturbation, since the variance tends to zero as $h$ is reduced.

\subsection{Robustness to Rotation}
\begin{figure}[H]
  \centering
  \includegraphics[width=0.75\textwidth]{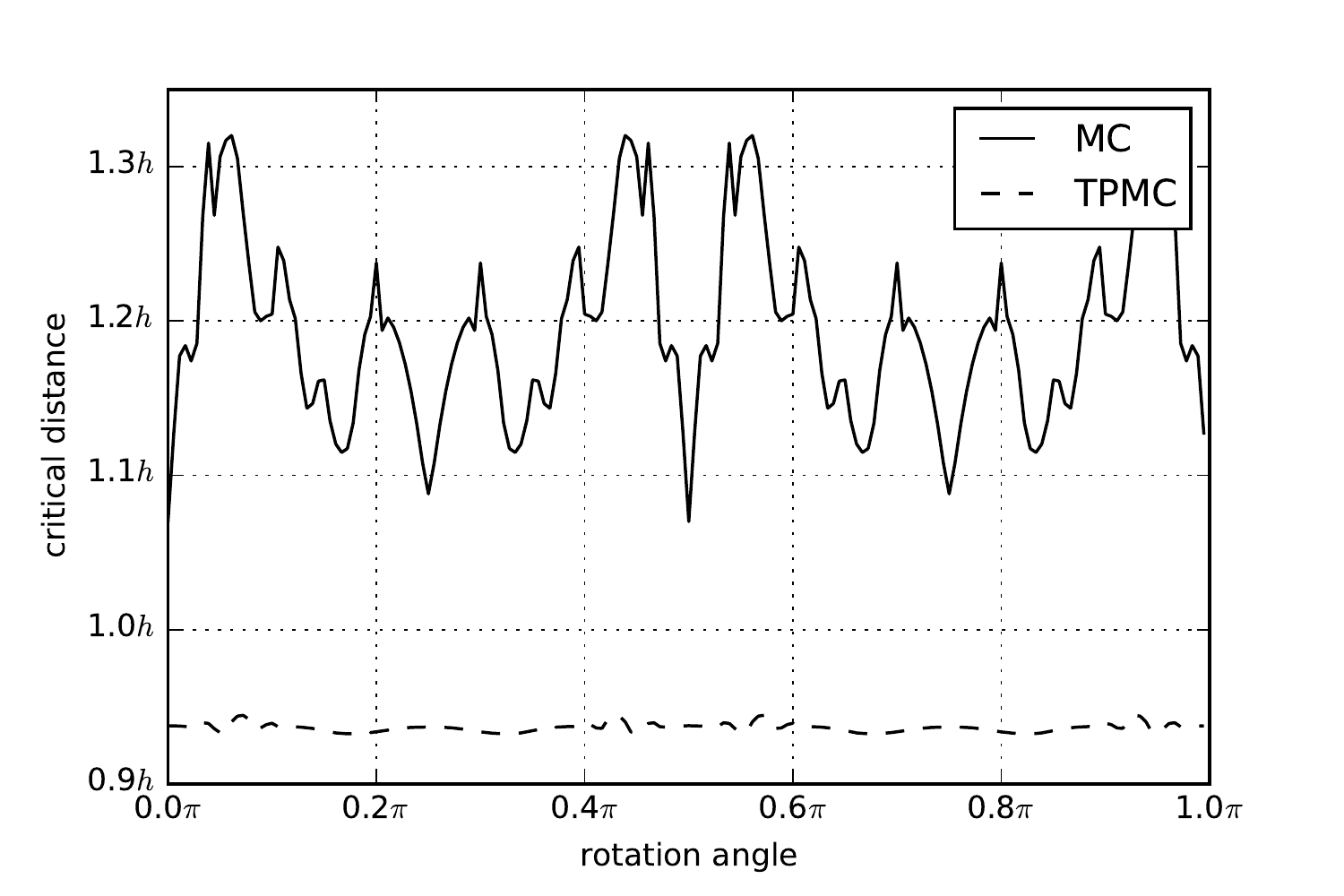}
  \caption{Two entangled tori are rotated in a three dimensional structured grid with rectangular elements. The distance between the tori is reduced and the distance at which the topology changes is noted in this figure as a multiple of $h$; smaller is better. The solid line represents the data obtained by \MC{} and the dashed one the data obtained by \ourMC{}. \MC{} is not robust with respect to the orientation of the tori.}
  \label{fig:plot_topology_change}
\end{figure}

For $d \rightarrow 0$ the topology changes, as the two tori ``merge''. In the discrete
reconstruction this happens already for $d > 0$.
In this test we compute, for different
directions $\phi$ of the two tori, the critical distance $d_c$ at which the
topology changes, i.e. the number of
discrete connected components switches from three to two.
In this test the grid size is fixed $h=\frac1{64}$ and $R=0.25$, as in
the first test.
For values of $\phi \in [0,\pi)$ we modify the
minor radius $r = \frac{R-d}{2}$ to compute the critical
distance $d_c = x\cdot h$, using a
bisection method.
Figure \ref{fig:plot_topology_change} shows the critical distance with (\ourMC) and
without (\MC) ambiguity-resolution (smaller is better).


For the method with ambiguity-resolution, the topology change occurs
relatively stable at a distance below $h$, yielding a sub-voxel resolution.
Without the ambiguity-resolution, meaning an anisotropic choice in an ambiguous case, the distance depends on the angle $\phi$ and is overall above $h$.

\subsection{Performance}

In order to evaluate the performance overhead of \ourMC{} as compared to \MC{}, we generate pseudo-random values in $[-1.0,1.0]$ for each vertex of a three dimensional structured grid.
We construct $30$ data sets for each grid size and compute the unique key for each grid cell using \ourMC{} and \MC{}.
For each data set, we average the time consumption over $30$ iterations.
We measure both the isolated time used for ambiguity resolution, and the total time consumption including the grid iteration.

This is a worst case scenario, as $\Phi$ is usually a smooth function,
so this test will encounter significantly more ambiguous case than
an average application.

The resulting ratios between time for \ourMC{} and \MC{} can be seen in Table \ref{tab:plot_time_comparison}.
The table also shows the number of face and center tests in relation to the number of key computations.
\begin{table}
  \centering
    \begin{tabular}{@{\quad}rr@{ $\pm$ }lr@{ $\pm$ }l@{\qquad}cc@{\quad}}
      \hline\hline
      \rule{0pt}{1.9ex}grid size
          & \multicolumn{2}{c}{\quad ambiguity resolution\qquad}
          & \multicolumn{2}{c}{total}
          & \multicolumn{2}{c}{test fraction}\\
          & \quad \phantom{ambi}$\overline{T}_\text{rel}$ & $\sigma_1$
          & $\overline{T}_\text{rel}$ & $\sigma_1$
          & face & center\\\hline
        16  & 2.27 & 2.07$\%$ & 1.23 & 2.57$\%$ & 76.0$\%$ & 22.0$\%$\\
        32  & 2.20 & 1.24$\%$ & 1.18 & 0.48$\%$ & 74.8$\%$ & 21.8$\%$\\
        64  & 2.17 & 0.72$\%$ & 1.17 & 0.25$\%$ & 74.9$\%$ & 21.8$\%$\\
        128 & 2.18 & 0.64$\%$ & 1.17 & 0.34$\%$ & 74.8$\%$ & 21.8$\%$\\
        256 & 2.19 & 1.86$\%$ & 1.18 & 0.50$\%$ & 74.8$\%$ & 21.8$\%$\\\hline\hline
    \end{tabular}
  \caption{Worst case test: Relative time consumption $\Tr=T_{\text{\ourMC{}}}/T_{\text{\MC{}}}$ for \ourMC{} as compared to
    \MC{} on a three dimensional structured grid with $N$ rectangular
    cells in each dimension. The timings are either the time used for
    ambiguity resolution or the total time, including grid
    iteration. We also show the fraction of tests involving face or center
    tests. For each grid size we compute 30 random datasets.}
  \label{tab:plot_time_comparison}
\end{table}
We observe an increased time consumption of \ourMC{} by a factor of approximately $2.2$ for the ambiguity resolution, which can be expected due to the need for additional tests.
For approximately three out of four key computations we need a face test and for two out of five we need a center test.
When taking also the iteration over the simple structured grid into
account, the factor already reduces to approximately $1.18$.
We note that for a computationally more expensive grid iteration, for example using an unstructured grid, the factor will further decrease.
As mentioned before, we want to point out that for smoother data, the need for
additional tests is in general significantly smaller than for random
data and therefore the performance overhead of the proposed method
will significantly smaller.

\section{Conclusion}
\label{sec:conclusion}

We proposed a method for geometric integration over implicitly described domains.
The method shows second order accuracy and a topologically correct representation of the geometry.
The method is robust with respect to grid perturbation while imposing only a small performance overhead when compared to a non-correct algorithm.
For an extension of the method to multiple level sets, one might consider a recursive application of the presented method.
For such a purpose, the interpolation of the vertices has to be modified to take higher order edges into account.

\bibliographystyle{plainurl}
\bibliography{literature}




\end{document}